\documentclass{article}
\usepackage{amsfonts}
\usepackage{graphicx}
\usepackage{bm}
\usepackage{amssymb}
\usepackage{amsmath}

\begin{document}

\begin{center}
{\Large \textbf{General Algorithmic Search} }

\bigskip

Sergio Hern\'{a}ndez$^{1}$, Guillem Duran$^{1}$, Jos\'{e} M. Amig\'{o}$^{2}$

\medskip

$^1$HCSoft Programaci\'{o}n, S.L., 30007 Murcia, Spain \\
$^2$Centro de Investigaci\'{o}n Operativa, Universidad Miguel Hern\'{a}ndez,\\ 
03202 Elche, Spain \\
E-Mails: sergio@hcsoft.net, guillem.db@gmail.com, jm.amigo@umh.es

\medskip
\end{center}

\begin{abstract}
In this paper we present a metaheuristic for global optimization called General Algorithmic Search (GAS). 
Specifically, GAS is a stochastic, single-objective method that evolves a swarm of agents in search of a global extremum. Numerical simulations with a sample of 31 test functions show that GAS outperforms Basin Hopping, Cuckoo Search, and Differential Evolution, especially in concurrent optimization, i.e., when several runs with different initial settings are executed and the first best wins. Python codes of 
all algorithms and complementary information are available online.

\medskip


\end{abstract}

\section{Introduction}
\label{intro}

This paper deals with unconstrained function optimization, i.e., with the
search of global maxima or minima of a single real-valued function without
constrains. General Algorithmic Search (GAS), the global optimization method
presented here, is a metaheuristic algorithm inspired by on-going work on
collective intelligence. Specifically, GAS is a stochastic, single-objective
method that evolves a swarm of agents in search of a global extremum. Its
performance is supported by numerical evidence and benchmarking against some
of the most popular metaheuristics currently used. For the present work we
have chosen three of them: Basin Hopping \cite{Wales1997}, Cuckoo Search 
\cite{Yang2009}, and Differential Evolution \cite{Storn1997}. For a review
of metaheuristics, see e.g. \cite{Voss2001}.

Optimization is a basic task in many areas of science, economy, and
technology, not least because many questions can be (re)formulated (usually
in a more elegant way) as optimization problems. As a well-known example,
the equations of motion of a mechanical system amount to the minimization of
its action --the so-called Principle of Least Action. In turn, some
optimization algorithms have been inspired by the traditional sciences,
e.g., Simulated Annealing \cite{Kirkpatrick1983}, Genetic Algorithms \cite%
{Holland1975}, Evolutionary Algorithms, and more. Nowadays optimization
remains a vigorous research field in which novel techniques are developed to
deal with old and new challenges such as those that arise in operations
research and machine learning. As the above few examples show, the
development of metaheuristics belongs certainly to this effort.

The rest of this paper is organized as follows. In Sec. 2 we describe the
GAS\ algorithm via a flowchart and comment some of its features. The
benchmark of GAS against Differential Evolution (DE), Cuckoo Search (CS),
and Basin Hopping (BH) is detailed in Sect. 3. The test-bed consists of 15
functions of two variables (see Annex) together with 16 higher-dimensional
functions, all of them deemed specially hard for the sake of global
optimization. The results of the benchmark are discussed in Sect. 4. They
show that GAS outperforms BH clearly, and CS and DE marginally in the case
of only one run. More importantly, the computational edge increases
significantly with the number of runs. We also zoom in on the optimization
of Lennard-Jones potentials because BH was specially developed to cope with
these functions. Here again the performance of GAS scales nicely with the
number of runs. The main messages of this paper are summarized in the
Conclusion (Sec. 5).

Finally, all the materials referred to in this paper are available online at
the URL \cite{Website}. In particular, Python codes of the above algorithms
as used in the computations reported here can be found in the file
Solvers.ipynb.

\section{The algorithm GAS}
\label{sec:1}

Given a continuous map $f:\mathbb{R}^{d}\rightarrow \mathbb{R}$, a typical
(single-objective) optimization problem has the form (without restriction)%
\begin{equation*}
\text{minimize }f(\mathbf{x})
\end{equation*}%
subject, in general, to some constraints on $\mathbf{x}%
=(x^{(1)},...,x^{(d)}) $. Thus, our objective is to find the position $%
\mathbf{x}^{\ast }$ of a global minimum of $f$ on a finite \textit{search
domain} $\mathcal{D}\subset \mathbb{R}^{n}$, i.e.,%
\begin{equation*}
\mathbf{x}^{\ast }=\arg \min_{\mathbf{x}\in \mathcal{D}}f(\mathbf{x})\text{,}
\end{equation*}%
where $\mathcal{D}$ might have a nontrivial topology due to hypothetical
constraints among the variables $x^{(n)}$, $1\leq n\leq d$. As a general
topological requisite we assume that $\mathcal{D}$ is compact (i.e., closed
and bounded), which guarantees the existence of $\mathbf{x}^{\ast }$. In
applications, $\mathcal{D}$ is usually an interval.

Roughly speaking, the GAS algorithm evolves a swarm of agents in pseudo-time
according to some rules. These agents will be called hereafter \textit{%
walkers} and they are characterized by both external and internal variables.
Specifically, the state of a walker comprises the current space coordinates
as well as the current value of its so-called \textit{flow}. Following the
basic idea of Tabu Search \cite{Glover1989,Glover1990}, local minima found
during the search will be flagged as `tabu' to avoid confinement in
suboptimal regions. As a new ingredient, the diffusion of the walkers in the
search domain is enhanced by a process called \textit{cloning} in which the
state of a given walker can be replaced by the state of another one.

Thus, the state of a walker $i$, $1\leq i\leq N$, is given by its position $%
\mathbf{x}_{i}=(x_{i}^{(1)},...,x_{i}^{(d)})\in \mathcal{D}$ and its flow $%
F_{i}>0$, to be defined in Eqn. (\ref{F_i}) below. For notational
convenience we denote the walker $i$ just by $\mathbf{x}_{i}$. Given a swarm
of $N$ walkers $\mathbf{x}_{i}$, set%
\begin{equation*}
f_{\min }=\min_{1\leq i\leq N}f(\mathbf{x}_{i}),\;\;\;f_{\max }=\max_{1\leq
i\leq N}f(\mathbf{x}_{i}),
\end{equation*}%
and%
\begin{equation}
\phi _{i}=\frac{f(\mathbf{x}_{i})-f_{\min }}{f_{\max }-f_{\min }}
\label{phi_i}
\end{equation}%
so that $0\leq \phi _{i}\leq 1$ for $i=1,...,N$, i.e., $\phi _{i}$ is a
scaling of $f(\mathbf{x}_{i})$. For the distance between walkers, we use the
Euclidean distance in $\mathbb{R}^{d}$:%
\begin{equation*}
dist(\mathbf{x}_{i},\mathbf{x}_{j})^{2}=\sum\nolimits_{n=1}^{d}\left(
(x_{i}^{(n)})^{2}-(x_{j}^{(n)})^{2}\right) .
\end{equation*}

The following flowchart describes how GAS searches for a global minimum $%
\mathbf{x}^{\ast }$ of $f$; some steps will be commented.

\begin{description}
\item[\textbf{Step 0. Initialization}] 

\item[0.1] The initial positions of the walkers, $\mathbf{x}_{i}$ ($1\leq
i\leq N$), are chosen randomly in the search domain $\mathcal{D}$.

\item[0.2] Compute%
\begin{equation*}
\mathbf{x}_{\min }=\arg \min_{1\leq i\leq N}\phi (\mathbf{x}_{i}).
\end{equation*}%
In the rare event of multiplicity, choose one of the minima randomly.

\item[0.3] Use $\mathbf{x}_{\min }$ and the algorithm L-BFGS-B \cite%
{Byrd1995} to find a local minimum $\mathbf{x}^{\ast }$.

\item[0.4] Fill out the \textit{tabu memory list} with $\mathbf{x}^{\ast }$,
i.e., set%
\begin{equation*}
\mathbf{t}_{r}=\mathbf{x}^{\ast },\;1\leq r\leq N\text{.}
\end{equation*}

\item[0.5] Set {\scriptsize BEST}$=\mathbf{x}^{\ast }$.
\end{description}

\textit{Comment}: The search loop starts at this point.

\medskip

\begin{description}
\item[Step 1. Walkers Flow \& Cloning] 

\item[1.1] For each walker $\mathbf{x}_{i}$ choose randomly another walker $%
\mathbf{x}_{j}$, $j\neq i$, and compute $d_{i,j}^{2}=dist(\mathbf{x}_{i},%
\mathbf{x}_{j})^{2}$.

\item[1.2] For each walker $\mathbf{x}_{i}$ choose randomly one tabu memory $%
\mathbf{t}_{r}$ and compute $\delta _{i,r}^{2}=dist(\mathbf{x}_{i},\mathbf{t}%
_{r})^{2}$. If $\mathbf{t}_{r}=\mathbf{x}_{i}$, set $\delta _{i,r}^{2}=1$.

\item[1.3] Compute the flow of $\mathbf{x}_{i}$, 
\begin{equation}
F_{i}=(\phi _{i}+1)^{2}\cdot d_{i,j}^{2}\cdot \delta _{i,r}^{2},  \label{F_i}
\end{equation}%
for $i=1,...,N$.

\item[1.4] For each walker $\mathbf{x}_{i}$ choose randomly another one $%
\mathbf{x}_{k}$, $k\neq i$, and compute the \textit{probability of} $\mathbf{%
x}_{i}$ \textit{cloning} $\mathbf{x}_{k}$:%
\begin{equation}
P_{i,k}=\left\{ 
\begin{array}{lc}
0 & \text{if }F_{k}>F_{i} \\ 
\min \{1,(F_{i}-F_{k})/F_{i}\} & \text{if }F_{k}\leq F_{i}%
\end{array}%
\right.   \label{P(i,k)}
\end{equation}%
Note that there is only one $k$ for each $i=1,...,N$.

\item[1.5] Take a random number $\rho \in \lbrack 0,1]$. If $\rho <$ $P_{i,k}
$ then the walker $\mathbf{x}_{i}$ copies the state of $\mathbf{x}_{k}$,
i.e., $\mathbf{x}_{i}\leftarrow \mathbf{x}_{k}$ and $F_{i}\leftarrow F_{k}$.
\end{description}

\bigskip

\textit{Comment}: Steps 1.3-1.5 promote the mobility of the walkers in such
a way that walkers can jump to distant, lower points of the landscape.

\medskip

\begin{description}
\item[Step 2. Local searches] 

\item[2.1] Find the `center of mass' of the walkers as follows:%
\begin{equation*}
\mathbf{x}_{cm}=\sum_{i=1}^{N}\phi _{i}\mathbf{x}_{i}.
\end{equation*}

\item[2.2] Compute 
\begin{equation*}
\mathbf{x}_{\min }=\arg \min_{1\leq i\leq N}\phi (\mathbf{x}_{i})
\end{equation*}%
(in the first loop, $\mathbf{x}_{\min }=\mathbf{x}^{\ast }$). In the rare
event of multiplicity, choose one of the minima randomly.

\item[2.3] Use the algorithm L-BFGS-B with $\mathbf{x}_{cm}$ and the current
swarm of walkers $\mathbf{x}_{i}$ to find a local minimum $\mathbf{t}$.

\item[2.4] Use the algorithm L-BFGS-B with $\mathbf{x}_{\min }$ and the
current swarm of walkers $\mathbf{x}_{i}$ to find another local minimum $%
\mathbf{t}^{\prime }$.

\item[2.5] Add $\mathbf{t}$ to the tabu memory list by overwriting a
randomly chosen tabu memory $\mathbf{t}_{i}$, and execute then the external
routine Memory Flow \& Cloning (see below).

\item[2.6] Add $\mathbf{t}^{\prime }$ to the tabu memory list by overwriting
a randomly chosen tabu memory $\mathbf{t}_{j}$, and execute then the
external routine Memory Flow \& Cloning (see below).

\item[2.7] Compute the minimum tabu memory, 
\begin{equation*}
\mathbf{x}^{\ast }=\min \{\mathbf{t}_{r}:1\leq r\leq N\},
\end{equation*}%
and update {\scriptsize BEST}$\leftarrow \mathbf{x}^{\ast }$.
\end{description}

\textit{Comment}: Steps 2.5 and 2.6 do not change the number of tabu
memories. These steps allow to find a global minimum which is close to a
local one.

\medskip

\begin{description}
\item[Halt criterion.] The usual criteria to exit a search loop include the
stability of the {\scriptsize BEST}, a maximum number of loops or function
reads, a maximum execution time, or possibly a mix of some of them. As
stability criterion we suggest to compare the current {\scriptsize BEST }%
with its average over a certain number of the preceding loops; exit then the
search loop at this point if the absolute value of the difference is smaller
than the precision sought. In any case, output%
\begin{equation*}
\mathbf{x}^{\ast }=\text{{\scriptsize BEST}}
\end{equation*}%
if the halt criterion is met; otherwise, continue.
\end{description}

\medskip

\begin{description}
\item[Step 3. Position update of the walkers and close the loop] 

\item[3.1] Define the `jump' of $\mathbf{x}_{i}$ as%
\begin{equation*}
\Delta _{i}=10^{-(5-4\phi _{i})}\in \lbrack 10^{-5},10^{-1}].
\end{equation*}

\item[3.2] If $L^{(n)}$ is the length of the $n$th dimension of $\mathcal{D}$
($1\leq n\leq d$), then update $\mathbf{x}_{i}$ as follows:%
\begin{equation*}
x_{i}^{(n)}\leftarrow x_{i}^{(n)}+L^{(n)}\xi ^{(n)},
\end{equation*}%
where each $\xi ^{(n)}$ is a random variable drawn from the normal
distribution $\mathcal{N}(0,\Delta _{i})$.

\item[3.3] If the resulting $\mathbf{x}_{i}$ does not belong to the search
domain $\mathcal{D}$, repeat Step 3.2 with $\Delta _{i}$ replaced by $\Delta
_{i}/2$ until $\mathbf{x}_{i}\in \mathcal{D}$.

\item[3.4] Go to Step 1.
\end{description}

\textit{Comment}: The scope of Steps 3.1-3.2 is that walkers with a small $%
\phi _{i}$ move slower than walkers with a large $\phi _{i}$, thus favoring
the accumulation of walkers on the lower zones of the landscape.

\textit{Comment}. In the case of searching for a global \textit{maximum }of $%
f$ (i.e., a global minimum of $-f$), the jump of $\mathbf{x}_{i}$ in Step
3.1 is defined by%
\begin{equation*}
\Delta _{i}=10^{-(1+4\phi _{i})}.
\end{equation*}

\medskip

\begin{description}
\item[External routine: Memory Flow \& Cloning] 

\item[R.1] For each tabu memory $\mathbf{t}_{r}$ choose randomly another
tabu memory $\mathbf{t}_{s}$, $s\neq r$, and compute $\tilde{d}%
_{r,s}^{2}=dist(\mathbf{t}_{r},\mathbf{t}_{s})^{2}$.

\item[R.2] Compute the flow of $\mathbf{t}_{r}$, 
\begin{equation}
\tilde{F}_{r}=(\tilde{\phi}_{r}+1)^{2}\cdot \tilde{d}_{r,s}^{2},  \label{F_r}
\end{equation}%
for $r=1,...,N$, where%
\begin{equation}
\tilde{\phi}_{r}=\frac{f(\mathbf{t}_{r})-\tilde{f}_{\min }}{\tilde{f}_{\max
}-\tilde{f}_{\min }},  \label{phi_r}
\end{equation}%
and 
\begin{equation*}
\tilde{f}_{\min }=\min_{1\leq r\leq N}f(\mathbf{t}_{r}),\;\;\;\tilde{f}%
_{\max }=\max_{1\leq r\leq N}f(\mathbf{t}_{r}).
\end{equation*}

\item[R.3] Compute the probability $\tilde{P}_{r,s}$ of cloning $\mathbf{t}%
_{s}$. To compute $\tilde{P}_{r,s}$ use formula (\ref{P(i,k)}) with $%
F_{i},F_{k}$ replaced by $\tilde{F}_{r},\tilde{F}_{s}$ respectively.

\item[R.4] Take a random number $\rho \in \lbrack 0,1]$. If $\rho <$ $\tilde{%
P}_{r,s}$ then $\mathbf{t}_{r}$ copies the state of $\mathbf{t}_{s}$, i.e., $%
\mathbf{t}_{r}\leftarrow \mathbf{t}_{s}$ and $\tilde{F}_{r}\leftarrow \tilde{%
F}_{s}$.
\end{description}

\textit{Comment}: Steps R.2-R.4 are meant to prevent static tabu memories
via cloning.

\medskip

\section{Benchmark}
\label{sec:2}

In this section, the performance of the algorithm presented in Sec. 2 is
compared against Basin Hopping (BH), Cuckoo Search (CS), and Differential
Evolution (DE). The reason for selecting BH and DE is that they are standard
algorithms of the Python libraries, what is a good indication of their
widespread use. As for CS, it was shown in (\cite{Yang2010}) to be
potentially more efficient that Swarm Particle Optimization \cite%
{Kennedy1995} and Genetic Algorithms. For the sake of this benchmark, we
boosted even more the efficiency of CS by adding to its standard
implementation a local search every 100 search loops. As mentioned in the
Introduction, the interested reader will find the Python codes used on the
website \cite{Website}\textbf{.}

The set of test functions is composed of 15 functions of two variables and
16 functions of more than two variables. Specifically, the 2-dimensional
test functions,%
\begin{equation}
f_{n}(x_{1},x_{2}),\;1\leq n\leq 15,  \label{2D test functions}
\end{equation}%
are the following: \textit{Ackley} ($n=1$), \textit{Beale} ($n=2$), \textit{%
Booth} ($n=3$), \textit{Easom} ($n=4$), \textit{Eggholder} ($n=5$), \textit{%
Goldstein-Price} ($n=6$), \textit{Levy N. 13} ($n=7$), \textit{Matyas} ($n=8$%
), \textit{MacCormick} ($n=9$), \textit{Rastrigin 2D} ($n=10$), \textit{%
Rosenbrock} ($n=11$), \textit{Schaffer N. 2} ($n=12$), \textit{Schaffer N. 4}
($n=13$), \textit{Sphere} ($n=14$), and \textit{Three-hump-camel} ($n=15$).
See the Annex for their definitions, search domains, and global minima.

For the higher dimensional test functions we use the following.

\begin{description}
\item[(i)] The 8 \textit{Lennard-Jones potentials} of $m$ particles, $3\leq
m\leq 10$, at positions $\mathbf{r}_{1},...,\mathbf{r}_{m}\in \mathbb{R}^{3}$
in reduced units (i.e., such that both the depth of the pair potential well
and the finite distance at which the pair potential vanishes are set equal
to $1$), 
\begin{equation}
L_{m}(\mathbf{r}_{1},...,\mathbf{r}_{m})=4\sum\limits_{1\leq i<j\leq
m}\left( \left( \frac{1}{r_{i,j}}\right) ^{12}-\left( \frac{1}{r_{i,j}}%
\right) ^{6}\right) ,  \label{L-J potentials}
\end{equation}%
where $r_{i,j}=dist(\mathbf{r}_{i},\mathbf{r}_{j})=r_{j,i}$. The number of
terms on the right hand side of (\ref{L-J potentials}) is $\frac{1}{2}m(m-1)$%
. The global minima of $L_{m}$ for $2\leq m\leq 110$ can be found in Table I
of \cite{Wales1997}. Here $\mathcal{D}=[-1.1,1.1]^{3m}$.

\item[(ii)] The 8 $d$-dimensional \textit{Rastrigin functions},%
\begin{equation*}
R_{d}(\mathbf{x})=10d+\sum_{n=1}^{d}\left( x_{n}^{2}-10\cos 2\pi
x_{n}\right) ,\;3\leq d\leq 10,
\end{equation*}%
with $\mathcal{D}=[-5.12,5.12]^{d}$, $\mathbf{x}^{\ast }=(0,...,0)$, and $%
R_{d}(\mathbf{x}^{\ast })=0$. As a matter of fact, $%
f_{10}(x_{1},x_{2})=R_{2}(x_{1},x_{2})$.
\end{description}

All the above test functions belong to a collection of functions generally
used for benchmarking in global optimization because of their challenging
geometrical properties; visit the website \cite{Website} for graphical
information. Furthermore, the Lennard-Jones potentials were included because
BH was originally designed to cope with these functions and shown in \cite%
{Wales1997} to improve the efficiency of genetic algorithms, the most
successful global optimization method for Lennard-Jones clusters till then.

As for the computational specifics, the calculations were done on an Intel
i7 with 8 cores, at 3.8 GHz, 16 GB of RAM, running Ubuntu 16.4. The halt
criterion was an error (absolute value of the difference between the current 
{\scriptsize BEST }and the known global minimum) $\leq 10^{-6}$.

\section{Results}
\label{sec:3}

Figure 1 summarizes the benchmark of GAS against BH, CS, and DE using the 31
multi-dimensional functions listed above as test functions; see the inlets
and caption for the color and line style codification.

First of all, panel (a) shows the percentage of global minima found with
each algorithm against the number of function reads. In each case, the code
was executed only once to generate the plots. In view of panel (a), we
conclude that GAS outperforms BH clearly and DE marginally, with CS
performing somewhere in-between.

More importantly, this favorable performance of GAS improves significantly
in the case of concurrent optimization. By these we mean that each algorithm
is run several times, each time with different initial settings (e.g.,
walker positions and random generator's seeds in the case of GAS), till a
given number of function reads has been achieved. Panels (b)-(d) plot the
highest percentage of minima found with $T=10,20,50$\ runs, respectively, vs
the number of function reads. Alternatively, one can think that each
algorithm is run simultaneously with $T$\ different initial settings. We see
at panels (b)-(d) that the GAS performance scales nicely with the number of
runs, while the CS and DE performances also improve but to a lesser extent
than the GAS performance does, and the BH performance hardly changes.

\begin{figure*}
  \includegraphics[width=1.0\textwidth]{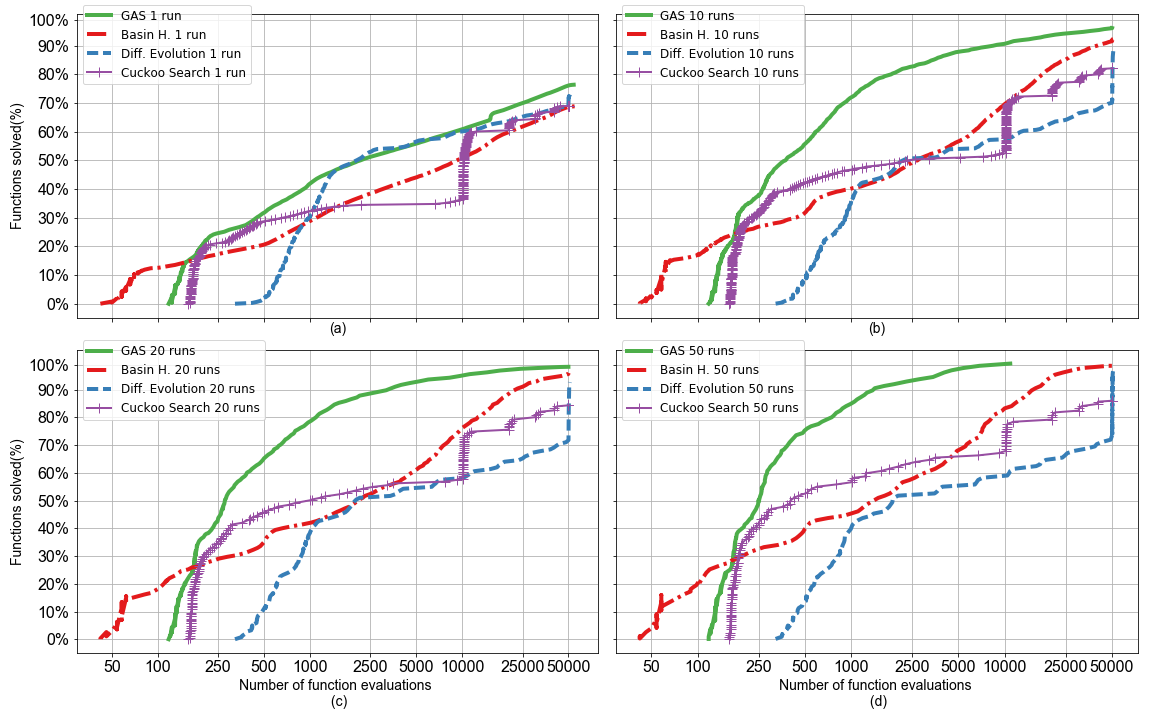}
\caption{(Color online) To compare the
performances of GAS (continuous thick line), BH (discontinuous segment-dot
thick line), CS (continuous thin line with transversal segments) and DE
(discontinuous segment-segment thick line) we plot the percentage of
optimized test functions against the number of function reads when executing
the algorithm one time (panel a). The same information is shown for
concurrent optimization with 10 runs (panel b), 20 runs (panel c) and 50
runs (panel d).}
\label{fig:1}       
\end{figure*}

As a complement to the above results, Figure 2 shows the performances of
GAS, BH, CS and DE when optimizing the Lennard-Jones potentials (%
\ref{L-J potentials}) of $m=4$, $6$, $8$, and $10$ particles at randomly chosen
positions (panels (a)-(d), respectively) with a single run. For the reason
mentioned in Sec. 3, we expected BH to score highest in this particular
benchmark but this was not the case. For $m=4$ the best performer is CS,
followed by GAS and BH. For $m=6$ GAS overtakes CS while BH remains third;
the performance gap between GAS and the other algorithms widens. Finally,
for $m=8,10$ BH overtakes CS but it still trails GAS at a considerable
distance; in both cases, the distance increases with the number of reads. DE
fails badly in all cases. Let us finally remark that GAS is an unbiased
algorithm, meaning that it does not exploit any specific insight into this
particular optimization problem nor into any other problem, for that matter.

\begin{figure*}
  \includegraphics[width=1.0\textwidth]{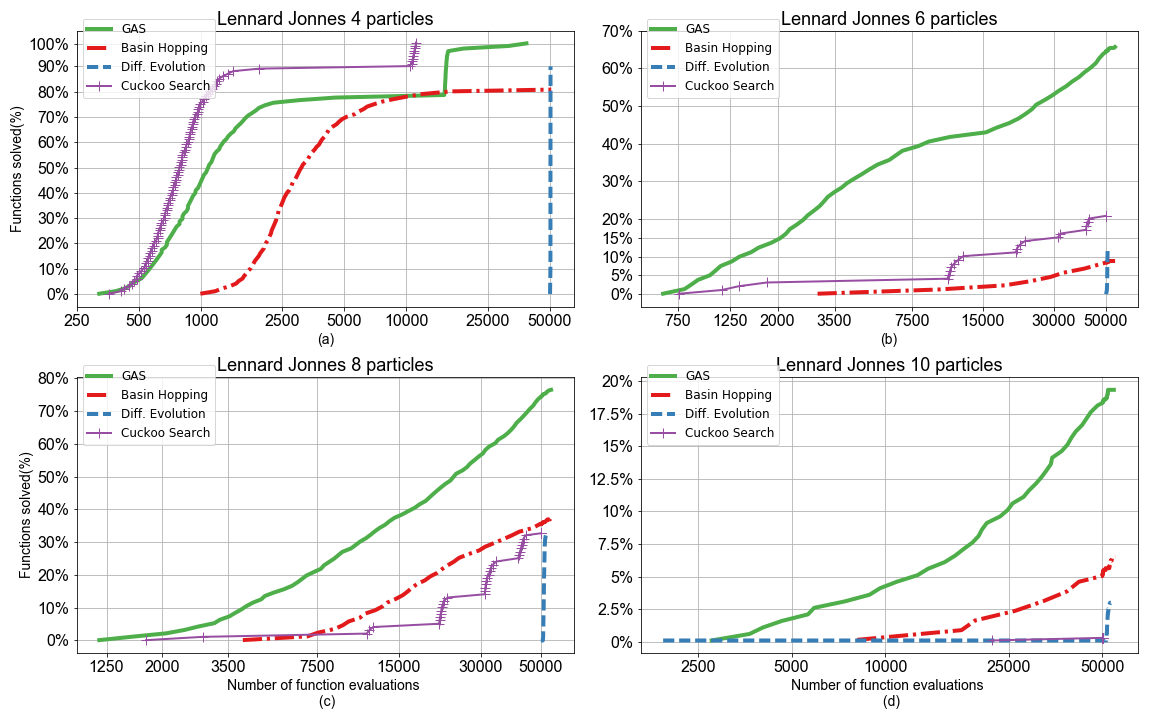}
\caption{(Color online) Line style
codification as in Fig. 1. Percentage of test functions optimized in one
single run against the number of function reads for the Lennard-Jones
potentials of 4 particles (panel a), 6 particles (panel b), 8 particles
(panel c), and 10 particles (panel d).}
\label{fig:2}       
\end{figure*}

\section{Conclusion}
\label{sec:4}

In this paper we presented a general metaheuristic for the search of global
extrema (GAS, Sec. 2). Among its features we highlight the use of walkers
with an internal label (called flow), randomization, and cloning (i.e., the
replacement of a walker or a tabu memory by a copy of another one).
Numerical simulations with a representative sample of test functions (Sect.
3) shows that GAS compares favorably to Basin Hopping, Cuckoo Search, and
Differential Evolution, the performance gap increasing with the number of
runs in a concurrent optimization (Sec. 4). We may conclude, therefore, that
GAS is a powerful tool for global optimization. The Python codes of the
algorithms used in the present work as well as additional information are
posted on the website \cite{Website}.

\bigskip

\begin{center}
\textbf{ANNEX: 2D test functions}

\medskip
\end{center}

Search domains ($\mathcal{D}\subset \mathbb{R}^{2}$), global minima ($%
\mathbf{x}^{\ast }\in \mathcal{D}$), and minimum values of the 2D test
functions (\ref{2D test functions}).

\begin{enumerate}
\item \textit{Ackley function}:%
\begin{equation*}
\begin{array}{c}
f_{1}(x_{1},x_{2})=-20\exp \left( -0.2\sqrt{0.5(x_{1}^{2}+x_{2}^{2})}\right)
\\ 
-\exp \left( 0.5(\cos 2\pi x_{1}+\cos 2\pi x_{2})\right) +e+20%
\end{array}%
\end{equation*}%
with $\mathcal{D}=[-5,5]\times \lbrack -5,5]$, $\mathbf{x}^{\ast }=(0,0)$,
and $f_{1}(\mathbf{x}^{\ast })=0$.

\item \textit{Beale function}:%
\begin{eqnarray*}
f_{2}(x_{1},x_{2})
&=&(1.5-x_{1}+x_{1}x_{2})^{2}+(2.25-x_{1}+x_{1}x_{2}^{2})^{2} \\
&&+(2.625-x_{1}+x_{1}x_{2}^{3})^{2}
\end{eqnarray*}%
with $\mathcal{D}=[-4.5,4.5]\times \lbrack -4.5,4.5]$, $\mathbf{x}^{\ast
}=(3,0.5)$, and $f_{2}(\mathbf{x}^{\ast })=0$.

\item \textit{Booth function}:%
\begin{equation*}
f_{3}(x_{1},x_{2})=(x_{1}+2x_{2}-7)^{2}+(2x_{1}+x_{2}-5)^{2}
\end{equation*}%
with $\mathcal{D}=[-10,10]\times \lbrack -10,10]$, $\mathbf{x}^{\ast }=(1,3)$%
, and $f_{3}(\mathbf{x}^{\ast })=0$.

\item \textit{Easom function}:%
\begin{equation*}
f_{4}(x_{1},x_{2})=-\frac{\cos x_{1}\cos x_{2}}{\exp \left( (x_{1}-\pi
)^{2}+(x_{2}-\pi )^{2}\right) }
\end{equation*}%
with $\mathcal{D}=[-100,100]\times \lbrack -100,100]$, $\mathbf{x}^{\ast
}=(\pi ,\pi )$, and $f_{4}(\mathbf{x}^{\ast })=-1$.

\item \textit{Eggholder function}:%
\begin{equation*}
\begin{array}{c}
f_{5}(x_{1},x_{2})=-(x_{2}+47)\sin \sqrt{\left\vert \frac{x_{1}}{2}%
+x_{2}+47\right\vert } \\ 
-x_{1}\sin \sqrt{\left\vert x_{1}-(x_{2}+47)\right\vert }%
\end{array}%
\end{equation*}%
with $\mathcal{D}=[-512,512]\times \lbrack -512,512]$, $\mathbf{x}^{\ast
}=(512,404.2319)$, and $f_{5}(\mathbf{x}^{\ast })=-959.6407$.

\item \textit{Goldstein-Price function}:%
\begin{eqnarray*}
&&f_{6}(x_{1},x_{2}) \\ 
&=&\left(1+(x_{1}+x_{2}+1)^{2}(19-14x_{1}+3x_{1}^{2}-14x_{2}+6x_{1}x_{2}+3x_{2}^{2})%
\right) \times \\
&&\times \left(
30+(2x_{1}-3x_{2})^{2}(18-32x_{1}+12x_{1}^{2}+48x_{2}-36x_{1}x_{2}+27x_{2}^{2})\right)
\end{eqnarray*}%
with $\mathcal{D}=[-2,2]\times \lbrack -2,2]$, $\mathbf{x}^{\ast }=(0,-1)$,
and $f_{6}(\mathbf{x}^{\ast })=3$.

\item \textit{Levy function N.13}:%
\begin{equation*}
\begin{array}{c}
f_{7}(x_{1},x_{2})=\sin ^{2}3\pi x_{1}+(x_{1}-1)^{2}(1+\sin ^{2}3\pi x_{2})
\\ 
+(x_{2}-1)^{2}(1+\sin ^{2}2\pi x_{2})%
\end{array}%
\end{equation*}%
with\ $\mathcal{D}=[-10,10]\times \lbrack -10,10]$, $\mathbf{x}^{\ast
}=(1,1) $, and $f_{7}(\mathbf{x}^{\ast })=0$.

\item \textit{Matyas function}:%
\begin{equation*}
f_{8}(x_{1},x_{2})=0.26(x_{1}^{2}+x_{2}^{2})-0.48x_{1}x_{2}
\end{equation*}%
with $\mathcal{D}=[-10,10]\times \lbrack -10,10]$, $\mathbf{x}^{\ast }=(0,0)$%
, and $f_{8}(\mathbf{x}^{\ast })=0$.

\item \textit{McCormick function}:%
\begin{equation*}
f_{9}(x_{1},x_{2})=\sin (x_{1}+x_{2})+(x_{1}-x_{2})^{2}-1.5x_{1}+2.5x_{2}+1
\end{equation*}%
with $\mathcal{D}=[-1.5,4]\times \lbrack -3,4]$, $\mathbf{x}^{\ast
}=(-0.54719,-1.54719)$, and $f_{9}(\mathbf{x}^{\ast })=-1.9133$.

\item \textit{Rastrigin 2D function}:%
\begin{equation*}
f_{10}(x_{1},x_{2})=20+x_{1}^{2}-10\cos 2\pi x_{1}+x_{2}^{2}-10\cos 2\pi
x_{2}
\end{equation*}%
with $\mathcal{D}=[-5.12,5.12]\times \lbrack -5.12,5.12]$, $\mathbf{x}^{\ast
}=(0,0)$, and $f_{10}(\mathbf{x}^{\ast })=0$.

\item \textit{Rosenbrock 2D function}:%
\begin{equation*}
f_{11}(x_{1},x_{2})=100(x_{2}-x_{1}^{2})^{2}+(x_{1}-1)^{2}
\end{equation*}%
with $\mathcal{D}=\mathbb{R}\times \mathbb{R}$, $\mathbf{x}^{\ast }=(1,1)$,
and $f_{11}(\mathbf{x}^{\ast })=0$.

\item \textit{Schaffer function N. 2}:%
\begin{equation*}
f_{12}(x_{1},x_{2})=0.5+\frac{\sin ^{2}(x_{1}^{2}-x_{2}^{2})-0.5}{%
(1+0.001(x_{1}^{2}+x_{2}^{2}))^{2}}
\end{equation*}%
with $\mathcal{D}=[-100,100]\times \lbrack -100,100]$, $\mathbf{x}^{\ast
}=(0,0)$, and $f_{12}(\mathbf{x}^{\ast })=0$.

\item \textit{Schaffer function N. 4}:%
\begin{equation*}
f_{13}(x_{1},x_{2})=0.5+\frac{\cos ^{2}\left( \sin \left\vert
x_{1}^{2}-x_{2}^{2}\right\vert \right) -0.5}{%
(1+0.001(x_{1}^{2}+x_{2}^{2}))^{2}}
\end{equation*}%
with $\mathcal{D}=[-100,100]\times \lbrack -100,100]$, $\mathbf{x}^{\ast
}=(0,1.25313)$, and $f_{13}(\mathbf{x}^{\ast })=0.292579$.

\item \textit{Sphere function}:%
\begin{equation*}
f_{14}(x_{1},x_{2})=x_{1}^{2}+x_{2}^{2}
\end{equation*}%
with $\mathcal{D}=\mathbb{R}\times \mathbb{R},$ $\mathbf{x}^{\ast }=(0,0)$,
and $f_{14}(\mathbf{x}^{\ast })=0$.

\item \textit{Three-hump-camel function}:%
\begin{equation*}
f_{15}(x_{1},x_{2})=2x_{1}^{2}-1.05x_{1}^{4}+\frac{x_{1}^{6}}{6}%
+x_{1}x_{2}+x_{2}^{2}
\end{equation*}%
with $\mathcal{D}=[-5,5]\times \lbrack -5,5]$, $\mathbf{x}^{\ast }=(0,0)$,
and $f_{15}(\mathbf{x}^{\ast })=0.$
\end{enumerate}


\end{document}